\documentclass{article}
\usepackage[utf8x]{inputenc}
\usepackage{amssymb,latexsym,amsmath,amsxtra,amsthm,amsfonts,amscd,slashed}
\usepackage{parskip}
\usepackage{enumerate}
\usepackage{color}
\usepackage{epsfig}
\usepackage{url}

\theoremstyle{plain}
\newtheorem{theorem}{Theorem}[section]

\newtheorem{lemma}[theorem]{Lemma}

\theoremstyle{remark}

\theoremstyle{definition}

\makeatletter
\def\thm@space@setup{
  \thm@preskip=\parskip \thm@postskip=0pt
}
\makeatother

\newcommand{\lap}[1]{\Delta#1}

\newcommand{\ind}{\mathbf{1}}
\newcommand{\at}[2]{\left.#1\right|_{#2}}
\newcommand{\paren}[1]{\left(#1\right)}
\newcommand{\norm}[1]{\Vert#1\Vert}

\newcommand{\abs}[1]{\left|#1\right|}
\newcommand{\brac}[1]{\left[#1\right]}

\newcommand{\vecs}[1]{\bf{#1}}
\newcommand{\ul}{u_{\lambda}}
\newcommand{\ext}[1]{E_0\brac{#1}}
\newcommand{\DO}[1]{\frac{d}{d#1}}

\DeclareMathOperator{\dist}{dist}

\def\R{\mathbb{R}}

\title{A Note on Moser Iteration and the Large Coupling Limit}
\date{}
\author{Ikemefuna Agbanusi\thanks{Department of Mathematics and Computer Science, Colorado College, \texttt{iagbanusi@coloradocollege.edu}}}
\begin{document}

\maketitle

\begin{abstract}
We consider heat semigroups of the form $\exp(t(\lap - \lambda\ind_{\Omega_0}))$ on bounded domains. These singularly perturbed equations arise in certain models of diffusion limited chemical reactions. Using variants of Moser iteration, we show sub-exponential decay in the so-called large coupling limit, i.e. as $\lambda\nearrow\infty$, in compact subdomains of the ``obstacle", $\Omega_0$. 
\end{abstract}

\section{Introduction}
Iterative and multi-scale methods are ubiquitous in analysis and applied mathematics, especially in the theory of singular perturbations. The main purpose of this article is to showcase another application of these ideas. In particular, we apply {\sc Moser's} \cite{Moser:1964rw} well known iteration technique to a parabolic P.D.E with a large parameter. We stress that the underlying ideas are not new, but we think leveraging Moser's scheme to obtain quantitative estimates for singularly perturbed PDE might be of wider interest. We turn now to a brief description of this problem.

Let $\Omega\subset\R^m$, $m\geq3$, be a bounded, open, connected subset with smooth boundary $\Gamma$. We are given a compact inclusion, $\Omega_0\Subset\Omega$, with boundary $\Gamma_0$ with $ \Omega_1:=\Omega\backslash \overline{\Omega}_0$. Also given is the Schr\"odinger operator $A_\lambda := \lap - \lambda\ind_{\Omega_0}$ with Neumann boundary conditions on $\Gamma$. For reasonable functions $f$, the large coupling problem is to determine the limit and convergence rate of $f(A_\lambda)$ as the coupling parameter $\lambda\nearrow\infty$. 

The case $f(x)=e^{tx}$, $t>0$, describes a diffusion process which is reflected on $\Gamma$ and removed with some probability on entering $\Omega_0$. It is well known that the corresponding large coupling limit is $e^{tB}$ where $B$ is the realization of $\lap$ in $L^2(\Omega_1)$ with Neumann boundary conditions at $\Gamma$ and Dirichlet boundary conditions at $\Gamma_0$. Physically, this limit corresponds to an instantaneous reaction or removal of the particle upon entering $\Omega_0$.

Using energy estimates, in \cite{Agbanusi:2016cpde} we showed an algebraic convergence rate in the exterior domain $\Omega_1$. In contrast, this paper is concerned with the behavior of solutions in the reaction region, $\Omega_0$.
Our first result is a sub-exponential convergence rate on compact, strictly interior subregions, $V$, of $\Omega_0$. Such decay is plausible in light of the Feynman--Kac formula. The novelty here is a proof using purely PDE techniques.
%
In the statement that follows, $E_0[g]$ is the extension of $g\in L^2(\Omega_1)$ by zero into $\Omega_0$ while $\ul(t,x)=e^{tA_\lambda}E_0[g]$.
\begin{theorem}\label{thm:main_result}
Assume $E_0[g]\in H^1(\Omega)$, and $r>0$ is the injectivity radius of $\partial V$. Set $a=\min\{\dist(\partial V,\Gamma_0),r\}$. For any $0<\nu<1/2$, there is a $\lambda_0=\lambda_0(a,\nu)$ such that for $\lambda\geq\lambda_0$,
\[ \sup_{t\geq0}\norm{\ul}^2_{L^2(V)}\leq e^{-\lambda^\nu}\lambda^{-1}\norm{\nabla g}^2_{L^2(\Omega_1)}.\]
\end{theorem}

Several comments are in order. First, the proof given in \S 2 also gives a similar estimate with a space--time $L^2$ norm: i.e. $L^2(I\times V)$ in place of $L^\infty(\R_+;L^2(V))$ for suitable intervals $I$. The parameter $a$ in Theorem \ref{thm:main_result} roughly measures the size of the boundary layer. The proof suggests this layer is $\sim\lambda^{-\frac12}$. From the modeling viewpoint, the estimate in Theorem \ref{thm:main_result} is related to the survival probability---the probability of finding the diffusing particle unreacted in $V$. Since $E_0[g]=0$ in $\Omega_0$, we assume the particle initially starts outside $\Omega_0$. Note that when $g\in H^1(\Omega_1)$ vanishes on $\Gamma_0$, then $E_0[g]\in H^1(\Omega)$. We denote this class of functions by $H_0^1(\Omega_1):=\{g\in H^1(\Omega_1): \at{g}{\Gamma_0}=0\}$.

The other main result is a \emph{pointwise} bound on solutions inside $\Omega_0$. 

%

\begin{theorem}\label{thm:2nd_result}
Assume in addition that $g\geq0$, then for $\lambda\geq\lambda_0$ as above and any $s\geq0$,
\[\sup_{Q_{a\gamma}} \ul^2\leq Ce^{-\lambda^\nu}\lambda^{-1}(a\gamma)^{-(m+2)}\norm{\nabla g}^2_{L^2(\Omega_1)},\]
where $Q_{a\gamma}=\{(t,x):t\in(s,s+(a\gamma)^2),\,x\in V_{-\gamma a}\}$, $0<\gamma<1/2$, and $C=C(m,V,\Omega_0)$.
\end{theorem}

The assumption $g\geq0$ is not essential, but it simplifies the argument and is natural in view of the probabilistic interpretation of the equations. The sets $V_{\rho}$ are defined in \S 2 and used there in the proof of Theorem \ref{thm:main_result}.
 The proof of Theorem \ref{thm:2nd_result} is sketched in \S 3. Both arguments employ energy estimates to a sequence of domains converging to $V$ which is the hallmark of Moser's scheme.
%

\section{Variations on a Theme of Moser and Sub-Exponential Decay}
Throughout we assume that the boundary, $\partial V$, is at least $C^1$, orientable and the Sobolev--Poincare inequalities hold on $V$. We begin with the observation that $\ul(t,x)=e^{tA_\lambda}E_0[g]$ solves:
\begin{align} \label{eq:sobolevGenpureDoi}
\partial_t\ul &=  \lap{\ul(t,x)} - \lambda \, \ind_{\overline{\Omega}_0}(x) \ul(t,x); \quad (t,x)\in I\times \Omega,
\end{align}
where $I=(0,T)$ with $0<T<\infty$. None of the constants in the estimates depend on $T$, so we may adjust it as appropriate. The initial and boundary conditions are:
\begin{equation}\label{eq:sobolevpureDoiBC}
   \left.\begin{aligned}
    &\at{\ul}{t=0} = \ext{g}:=
     \begin{cases}
    g(x),  &x \in \Omega_1;\\
    0, & x \in \Omega_0.
  \end{cases}\\
    & \nabla\ul(t,x)\cdot \hat{\vecs{n}} = 0,\quad\quad (t,x)\in I\times \Gamma.
 \end{aligned}
 \right\}
 \end{equation}

We need two estimates which we state as Lemmas. The first is contained in \cite{Agbanusi:2016cpde}. 
\begin{lemma}\label{lem:base_est}
If $g\in H_0^1(\Omega_1)$, then 
\[\sup_{t\geq0}\norm{\ul(t,\cdot)}^2_{L^2(\Omega_0)}\leq\frac{1}{\lambda}\norm{\nabla g}^2_{L^2(\Omega_1)};\qquad \norm{\ul}^2_{L^2(Q_0)}\leq\frac{1}{2\lambda}\norm{g}^2_{L^2(\Omega_1)}\]
\end{lemma}
The proof is a standard energy argument: multiply the PDE by the solution or one of its derivatives---in this case $\ul$ and $\partial_t\ul$---and integrate in space, time or both. Similar arguments are shown later so we omit this proof.

The other ingredient is a Caccioppoli type estimate. These are usually given on cubes, balls and their parabolic counterparts. Our unusual statement requires some notation. Viewing $\partial V$ as a $C^1$ compact, embedded, oriented surface in $\R^m$, as in {\sc lee} \cite[pg. 255--257]{Lee:2000vn}, we can define the map 
\begin{align}
\exp^{\perp}:\partial V\times \R\to \R^m;\qquad (x,s)\mapsto x+s\hat{\vecs{n}}(x),
\end{align}
where $\hat{\vecs{n}}(x)$ is the unit outward normal vector field on $\partial V$. Let $\partial V_\rho=\exp^{\perp}(\partial V,\rho)$ and define 
\[V_\rho=
\begin{cases}
V\bigcup\paren{\bigcup\limits_{0\leq s<\rho}\exp^{\perp}(\partial V,s)}, & \rho>0;\\
V\backslash \paren{\bigcup\limits_{\rho\leq s<0}\exp^{\perp}(\partial V,s)}, & \rho<0;
\end{cases}.\]
The injectivity radius $r$ is the supremum over all $|\rho|>0$ for which $\partial V$ is $C^1$ diffeomorphic to $\partial V_\rho$. Note that $\dist(V_\rho,V_{\rho+\sigma})$ is  $|\sigma|$ as can be shown using Gauss' Lemma.

\begin{lemma}\label{lem:caccio}
With $a=\min\{\dist(\partial V,\Gamma_0),r\}$ and for any $0\leq\rho<\rho+\sigma <a,$
\begin{align*}
&\norm{\ul}^2_{L^{\infty}(\R_+;L^2(V_\rho))}\leq\frac{4}{\lambda\sigma^2} \norm{\ul}^2_{L^{\infty}(\R_+;L^2(V_{\rho+\sigma}))};\\
&\norm{\ul}^2_{L^2(I\times V_\rho)}\leq\frac{4}{\lambda\sigma^2} \norm{\ul}^2_{L^2(I\times V_{\rho+\sigma})}.
\end{align*}
\end{lemma}
Before proving this lemma, let us show how it implies the theorem.
\begin{proof}[{\bf Proof of Theorem~\ref{thm:main_result}}.] Define $U_0:=\Omega_0$ and $U_j = V_{\gamma a(1-j/N)}$ for $1\leq j\leq N$, any $0<\gamma<1$ and some integer $N$ to be determined soon. Thus
\[V=:U_N\Subset U_{N-1} \Subset \ldots\Subset U_1 \Subset U_0:=\Omega_0,\]
with $\dist(U_j,U_{j-1})\geq\frac{\gamma a}{N}$.
With $X(U)$ either of the spaces $L^\infty(\R_+;L^2(U))$ or $L^2(I\times L^2(U))$,  Lemma \ref{lem:caccio} implies
\[\norm{u}^2_{X(U_j)}\leq\frac{4}{\lambda(\dist(U_j,U_{j-1}))^2}\norm{u}^2_{X(U_{j-1})}.\]
Iterating this $N$ times and applying Lemma \ref{lem:base_est} gives
\[\norm{u}^2_{X(U_N)}\leq\paren{\frac{4}{\lambda(\gamma a/N)^2}}^N\norm{u}^2_{X(U_0)}\leq\frac{\norm{g}_{H^1(\Omega_1)}}{\lambda}\paren{\frac{4N^2}{\gamma^2a^2\lambda}}^N,\]
Choosing $N$ as the integer part of $\lambda^\nu$ with $0<\nu<1/2$,  we see that
$\paren{\frac{4N^2}{\gamma^2a^2\lambda}}^N\leq\paren{\frac{4\lambda^{2\nu}}{\gamma^2a^2\lambda}}^{\lambda^\nu}$.
 The theorem follows once we realize that $\frac{4\lambda^{2\nu}}{\gamma^2a^2\lambda}\leq e^{-1}$ for $\lambda\geq\paren{\frac{4e}{\gamma^2a^2}}^{\tfrac{1}{1-2\nu}}$.
 \end{proof}
 
We can fiddle with $N$ as long as $N(\lambda)=o(\sqrt{\lambda})$ at the possible cost of a larger $\lambda_0$. One amusing example is $N(\lambda)=\lambda^{\frac{1}{2}-\frac{1}{2k}}(\log\lambda)^{\frac{k}{2}}$ for any natural number $k$. With so much leeway, the best choice of $N(\lambda)$ is an interesting question.


Right now, we turn to verifying Lemma~\ref{lem:caccio}.
\begin{proof}[{\bf Proof of Lemma~\ref{lem:caccio}}.] 
Let $0\leq\eta(x)\leq1$  be a smooth function with $\eta=1$ on $V_{\rho}$ and vanishing outside $V_{\rho+\sigma}$ chosen so that $\abs{\nabla\eta}\leq2/\sigma$, we get 
We multiply the PDE by $\eta^2\ul$ and integrate by parts using the fact that $\eta^2\ul$ is compactly supported in $\Omega_0$ to get
\begin{align*}
\int_{\Omega_0}\partial_t\ul\ul\eta^2\,dx &= \int_{\Omega_0}\lap\ul\ul\eta^2\,dx -\lambda \int_{\Omega_0}\ul^2\eta^2\,dx,\\
&=-\int_{\Omega_0}\nabla\ul\cdot\nabla(\ul\eta^2)\,dx -\lambda \int_{\Omega_0}\ul^2\eta^2\,dx, \\
&=-\int_{\Omega_0}\abs{\nabla\ul}^2\eta^2\,dx-\int_{\Omega_0}2\eta\ul\nabla\eta\cdot\nabla\ul\,dx -\lambda \int_{\Omega_0}\ul^2\eta^2\,dx.
\end{align*}
For simplicity, we drop the integration measure. Cauchy's inequality ``with $\epsilon$" implies
\begin{align*}
\frac{1}{2}\DO{t}\int_{\Omega_0}\ul^2\eta^2 +\lambda\int_{\Omega_0}\ul^2\eta^2+\int_{\Omega_0}\abs{\nabla\ul}^2\eta^2&\leq2\abs{\int_{\Omega_0}\eta\ul\nabla\eta\cdot\nabla\ul},\\
&\leq \frac{1}{2}\int_{\Omega_0}\eta^2\abs{\nabla\ul}^2+2\int_{\Omega_0}\ul^2\abs{\nabla\eta}^2,
\end{align*}
which after some rearrangement gives
\begin{equation}\label{eqn:main_L2_estimate}
\frac{1}{2}\DO{t}\int_{\Omega_0}\ul^2\eta^2 +\lambda\int_{\Omega_0}\ul^2\eta^2+\frac{1}{2}\int_{\Omega_0}\abs{\nabla\ul}^2\eta^2\leq2\int_{\Omega_0}\ul^2\abs{\nabla\eta}^2.
\end{equation}
The positivity of the term involving the gradient on the left side implies that 
\[\frac{1}{2}\DO{t}\int_{\Omega_0}\ul^2\eta^2 +\lambda\int_{\Omega_0}\ul^2\eta^2\leq2\int_{\Omega_0}\ul^2\abs{\nabla\eta}^2,\]
and Gronwall's Lemma shows that
\[\norm{(\ul\eta)(t)}^2_{L^2(\Omega_0)}\leq   4\int_0^t\int_{\Omega_0}e^{-2\lambda(t-s)}\ul^2\abs{\nabla\eta}^2\,dx\,dt .\]
Using the fact that $\abs{\nabla\eta}\leq2/\sigma$, we get 

\[\norm{\ul(t)}^2_{L^2(V_\rho)}\leq\norm{(\ul\eta)(t)}^2_{L^2(\Omega_0)}\leq\frac{4}{\lambda\sigma^2}\norm{\ul}^2_{L^\infty(\R_+;L^2(V_{\rho+\sigma}))},\]
and taking the supremum over $t$ on the left proves the ``$L^\infty$ part" of the Lemma.

To get the ``$L^2$ part", we integrate \eqref{eqn:main_L2_estimate} from 0 to $T$ to get
\[\frac{1}{2}\norm{\ul(t)}^2_{L^2(V_\rho)}+\lambda\norm{\ul}^2_{L^2(I\times V_\rho)}+\frac{1}{2}\norm{\eta\nabla\ul}^2_{L^2(I\times V_\rho)}\leq \int_0^T\int_{\Omega_0}\ul^2\abs{\nabla\eta}^2\,dx\,dt,\]
from which we get
\[\norm{\ul}^2_{L^2(I\times V_\rho)}\leq\frac{4}{\lambda\sigma^2}\norm{\ul}^2_{L^2(I\times V_{\rho+\sigma})}.\]
\end{proof}

Though we squeezed two estimates from the proof above, it has a bit more to give:
\begin{enumerate}[(i)]
\item One refinement---which we won't pursue here---is to capitalize on the fact that $\nabla\eta$ is actually supported in $V_{\rho+\sigma}\backslash V_\rho$ and not on the whole of $V_{\rho+\sigma}$.
\item Using the $L^\infty$ estimate in $t$ after the application Grownwall's Lemma does not take advantage of the explicit integral representation.
\item We completely ignored the gradient term in the energy inequalities.
\end{enumerate}

The third issue is examined in \S3 with the help of the Poincare--Sobolev inequality which relates $L^p$ norms of the function and its gradient.

\section{Moser Iteration all the Way}
This section is devoted to the proof of Theorem \ref{thm:2nd_result}. It is implied by the following ``mean value inequality":
\[\sup_{Q_\gamma} \ul^2\leq C\iint\limits_{(s,s+a^2\gamma^2)\times V}\ul^2\,dx\,dt.\]
The argument is simple: the quantity on the right is controlled by $\norm{\ul}^2_{L^2(I\times V)}$ which decays sub-exponentially by Theorem \ref{thm:main_result}.
The inequality is so named  because in the case of a sphere, the constant, $C$ is proportional to the reciprocal of the volume of the region of integration.

The proof of the mean value inequality follows {\sc moser}'s proof of a similar result in \cite{Moser:1964rw} (see the exposition in {\sc saloff-coste} \cite[pg 445-447]{Saloff-Coste:1995uf}). Partly because our formulation is slightly different, and partly to keep this note self--contained, we include a sketch of the proof. 

We begin with the following identity for $p\geq1$ and non-negative $\ul$
\begin{multline*}\frac{1}{2p}\DO{t}\int\limits_{\Omega_0}(\ul^{p}\chi\eta)^2+\chi^2\int\limits_{\Omega_0}|\nabla(\ul^{p}\eta)|^2+\chi^2\paren{1-p}^2\int\limits\limits_{\Omega_0}\ul^{2p-2}\eta^2\abs{\nabla \ul}^2+\lambda\chi^2\int\limits_{\Omega_0}\ul^{2p}\eta^2\\
=\chi'\chi\int\limits_{\Omega_0}\ul^{2p}\eta^2+\chi^2\int\limits_{\Omega_0}|\nabla\eta|^2\ul^{2p},
\end{multline*}
where $\chi=\chi(t)$ is a smooth function of $t$ and $\eta$ we have encoutered before. This can be derived by multiplying the differential equation by $\ul^{2p-1}\chi^2\eta^2$ and integrating by parts. We omit the details, but when the algebraic dust settles the above identity emerges. Note that when $p=1$, $\chi\equiv1$, Cauchy's inequality allows us to recover the main inequality, \eqref{eqn:main_L2_estimate}, in the proof of Lemma \ref{lem:caccio}.

It's time to specify $\chi$ and $\eta$. With $0<\sigma<\widetilde{\sigma}$, we choose $0\leq\eta(x)\leq1$ with $\eta=1$ on $V_{\rho +\sigma}$ and vanishing outside $V_{\rho+\widetilde{\sigma}}$ and $0\leq\chi(t)\leq1$ with $\chi=1$ on $(-\infty,s+\sigma^2)$ and vanishing outside $(s+\widetilde{\sigma}^2,\infty)$. Let $I_\sigma=(s,s+\sigma^2)$ and integrate the identity from $s$ to any $t\in I_\sigma$ to get
\begin{equation}\label{eqn:main_p_estimate}
\sup_{t\in I_{\sigma}}\int\limits_{V_{\rho+\sigma}}\ul^{2p}+2p\iint\limits_{I_\sigma\times V_{\rho+\sigma}}|\nabla(\ul^{p})|^2\leq \frac{8p}{(\widetilde{\sigma}-\sigma)^2}\iint\limits_{I_{\widetilde{\sigma}}\times V_{\rho+\widetilde{\sigma}}}\ul^{2p}
\end{equation}
We have used the support properties of $\chi$ and $\eta$ and the fact that $|\chi'|\leq 2/(\widetilde{\sigma}-\sigma)^2$ and $|\nabla\eta|^2\leq 2/(\widetilde{\sigma}-\sigma)^2$. The $p$ in the fraction on the right side of the inequality can probably be eliminated, but we proceed as is. Later on, we will use $Q_{\sigma}:=I_\sigma\times V_{\rho+\sigma}$ to further simplify the notation.

H\"older's inequality implies that
\[\int_{U}w^{2(1+\frac{2}{m})} \leq\paren{\int_{U}w^{2}}^{\frac{2}{m}}\paren{\int_{U}w^{\frac{2m}{m-2}}}^{\frac{m-2}{m}},\]
while the Poincare--Sobolev inequality for compactly supported functions implies 
\[\paren{\int_{U}w^{\frac{2m}{m-2}}}^{\frac{m-2}{m}}\leq \kappa\int_{U}|\nabla w|^2,\]
and $\kappa=\kappa(U,m)$ is the best constant for the inequality. Hence, for any interval $J$,
\[\iint\limits_{J\times U}w^{2(1+\frac{2}{m})} \leq \kappa\sup_{t\in J}\paren{\int_{U}w^{2}}^{\frac{2}{m}}\iint\limits_{J\times U}|\nabla w|^2.\]
Applying the above inequality with $J=I_\sigma$ and $U=V_{\rho+\sigma}$  and using \eqref{eqn:main_p_estimate}
gives
\begin{equation}\label{eqn:recursion}
 \iint\limits_{Q_{\sigma}}\ul^{2p\theta}\leq2^{\frac{2}{m}}\kappa \paren{\frac{4p}{(\widetilde{\sigma}-\sigma)^2}\iint\limits_{Q_{\widetilde{\sigma}}}\ul^{2p}}^\theta,
 \end{equation}
with $\theta=1+\tfrac{2}{m}$ and $\kappa$ is the Poincare--Sobolev constant for $V$.
To iterate this, we put $\rho=-2a\gamma$, with $0<\gamma<1/2$, and define $\sigma_i=a\gamma(1+2^{-i})$ for $i=0,1,2\ldots$ Note that $\sigma_0=2a\gamma$ and $\lim\limits_i\sigma_i=a\gamma$. Applying \eqref{eqn:recursion} with $p=\theta^i$ and $\sigma=\sigma_{i+1}$ and $\widetilde{\sigma}=\sigma_i$ we get:
\[\iint\limits_{Q_{\sigma_{i+1}}}\ul^{2\theta^{i+1}}\leq2^{\frac{2}{m}}\kappa \paren{\frac{4\theta^i4^{i+1}}{(a\gamma)^2}\iint\limits_{Q_{\sigma_i}}\ul^{2\theta^i}}^\theta.\]
It then follows that
\[\paren{\,\,\iint\limits_{Q_{\sigma_{i+1}}}\ul^{2\theta^{i+1}}}^{\theta^{-i-1}}\leq (2^{\frac{2}{m}}\kappa)^{\sum_{j=1}^{i+1}\theta^{-j}}\paren{\frac{16}{(a\gamma)^2}}^{\sum_{j=0}^i\theta^{-j}}(4\theta)^{\sum_{j=1}^{i}j\theta^{-j}}\iint\limits_{Q_{\sigma_0}} \ul^2.\]
Sending $i\nearrow\infty$ and playing with several geometric series gives
\[\sup\limits_{Q_{a\gamma}}\ul^2\leq \frac{C(m)\kappa^{\frac{m}{2}}}{(a\gamma)^{m+2}}\iint\limits_{Q_{\sigma_0}}\ul^{2},\]
for some explicit constant $C(m)$ depending only on $m$. This is the desired inequality.

\subsection*{Acknowledgement} We would like to thank Professors Gene Wayne, Samuel Isaacson and Andres Larrain-Hubach for providing feedback on an earlier draft.


\thispagestyle{empty}
\bibliography{lib_papers} 
\bibliographystyle{amsplain}
\end{document}